\documentclass[twoside,leqno,twocolumn,letter]{article}  
\usepackage{ltexpprt,url,amsmath,amsfonts,amssymb,graphics,graphicx,color} 
\usepackage[margin=1.0in]{geometry}

\begin{document}

\title{\bf \Large On Time-optimal Trajectories for a Car-like Robot with One Trailer}
\author{Hamidreza Chitsaz\thanks{Department of Computer Science, Wayne State University, 5057 Woodward Ave, Suite 3010, Detroit, MI 48202.}\\ \url{chitsaz@wayne.edu}}
\date{}

\maketitle


\begin{abstract} \small\baselineskip=9pt 
In addition to the theoretical value of challenging optimal control problmes, recent progress in autonomous vehicles mandates further research in optimal motion planning for wheeled vehicles. Since current numerical optimal control techniques suffer from either the curse of dimensionality, e.g. the Hamilton-Jacobi-Bellman equation, or the curse of complexity, e.g. pseudospectral optimal control and max-plus methods, analytical characterization of geodesics for wheeled vehicles becomes important not only from a theoretical point of view but also from a practical one. Such an analytical characterization provides a fast motion planning algorithm that can be used in robust feedback loops. In this work, we use the Pontryagin Maximum Principle to characterize extremal trajectories, i.e. candidate geodesics, for a car-like robot with one trailer. We use time as the distance function. In spite of partial progress, this problem has remained open in the past two decades. Besides straight motion and turn with maximum allowed curvature, we identify planar elastica as the third piece of motion that occurs along our extremals. We give a detailed characterization of such curves, a special case of which, called \emph{merging curve}, connects maximum curvature turns to straight line segments. The structure of extremals in our case is revealed through analytical integration of the system and adjoint equations. 
\end{abstract}

\section{Introduction}
With the advent of autonomous vehicles, there is a renewed need for \emph{analytical} characterization of geodesics for wheeled vehicles, as this is an important problem not only from a theoretical point of view but also from a practical one. Such characterizations provide motion planning algorithms with constant time complexity that can be used in a robust feedback loop. Moreover, geodesics for wheeled vehicles provide a library of motion primitives from which collision-free (locally) optimal paths can be constructed in the presence of obstacles. Existing alternative methods, namely numerical optimal control techniques, suffer from either the curse of dimensionality, e.g. the Hamilton-Jacobi-Bellman equation requires $O(n^d)$ space for a $d$-dimensional system \cite{Ber87}, or the curse of complexity, e.g. pseudospectral optimal control algorithms that solve a general finite dimensional nonlinear optimization per initial-goal pair \cite{ElnKazRaz95,GarPatHagRaoBenHun10,RosFah04,RosKar12} and max-plus methods that approximate the solution of the Hamilton-Jacobi-Bellman equation, without an explicit representation of the discretized configuration space \cite{GauMceQu11,Mce05}. 

\begin{figure}
\begin{center}
\includegraphics[scale=0.52]{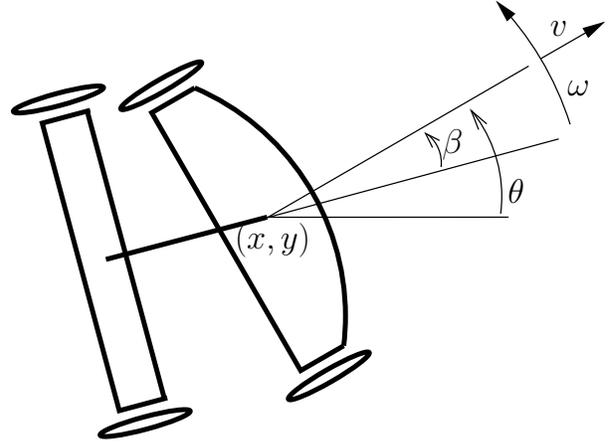}
\end{center}
\caption{The configuration $(x, y, \theta, \beta)$ and controls $(v, \omega)$ of a car-like robot with one trailer.} \label{fig:config}
\end{figure}

In particular, characterization of time optimal paths for a car-like robot with trailers as the premier type of geodesics has remained an open problem in spite of partial progress in the past two decades. Chyba and Sekhavat gave the first approach toward characterization of time optimal paths for a mobile robot with one trailer \cite{ChySek99}. They applied the Pontryagin Maximum Principle and switching structure equations as a necessary condition for optimality to give a partial characterization of extremals, i.e. those trajectories that satisfy the necessary condition and hence are candidate geodesics \cite{PonBolGamMis62,SusTan91}. In this paper, we continue that work to give a complete characterization of extremals with analytical solution for the system and adjoint differential equations. In particular, we demonstrate that a special maneuver, which we call \emph{merging}, comes in between a minimum radius turn and straight line motion along time optimal paths to straighten the trailer completely before merging into the 
straight line segment. It has not escaped our attention that the merging curve here happens to be the same as planar elastica \cite{ChiLav07b,Jur92,WalMonSas94}. This interesting coincidence requires further investigation.

\subsection{Related Work.}
In a classical paper in 1957, L.E. Dubins characterized the shortest paths, measured in the Euclidean plane, between two
points in an obstacle-free $SE(2)$, with the constraint that the average curvature of the projection of the path onto the Euclidean plane over any interval along the path
be bounded from above \cite{Dub57}. The Dubins
curves are composed of sequences of up to three primitives,
which consists in a line and an arc of a circle with maximum allowed curvature.
They are of particular interest in control and robotics because they give the shortest paths for a simple
airplane model as well as a forward-only car with bounded steering angle and velocity.
Extending the car model to allow backward motion as well, Reeds
and Shepp characterized the shortest paths for a car with
bounded steering angle and velocity \cite{ReeShe90}. Sussmann and Tang
developed a general methodology and machinery based on geometric optimal control 
for solving problems of
this sort by revisiting the Reeds-Shepp problem \cite{SusTan91}. Sou\`{e}res, Boissonnat, and Laumond gave an optimal control synthesis, i.e. a mapping from pairs of
initial-goal configurations to the optimal trajectories \cite{SouBoi98,SouLau96}. These developments lead in the new millenium to the characterization of time optimal trajectories
for the differential drive
\cite{BalMas02} and omni-directional vehicles in the plane \cite{BalKavMas08,FurBal10,FurBalChiKav08}. Chyba \emph{et al.} considered time optimal trajectories for an underwater vehicle \cite{ChyHab05,ChyHabSmiCho08}.
We characterized time optimal trajectories for a simple 3D airplane model obtained by extending the Dubins car model \cite{ChiLav07b} and geodesics for a differential
drive robot with the total wheel rotation as the distance function \cite{ChiLav07a,ChiLavBalMas09}.

\section{Problem Specification} In this paper, we consider the analytical solution of the following optimal control problem for every initial and goal configuration pair $q_i, q_g$:
\begin{equation}
\min \int_0^T dt
\end{equation}
subject to $q(0) = q_i$, $q(T) = q_g$, and
\begin{equation}\label{equ:system}
\dot q = \left( \begin{array}{c}
\dot x \\
\dot y \\
\dot \theta \\
\dot \beta
\end{array} \right) = \left( \begin{array}{c}
v \cos\theta \\
v \sin\theta \\
\omega \\
-v \sin\beta + \omega
\end{array}\right),
\end{equation}
in which $q=(x, y, \theta, \beta) \in \mathbb{R}^2\times S^1\times S^1$ is the configuration and $u=(v, \omega) \in U = [-1, 1]^2$ is the control. Figure \ref{fig:config} depicts the configuration of the robot, where $v$ is the linear velocity, $\omega$ is the angular velocity, $(x, y)$ is the coordinate of the center of robot, $\theta$ is the orientation, and $\beta$ is the angle between the trailer and robot.

\section{Pontryagin Maximum Principle}
Let $\lambda = (\lambda_x, \lambda_y, \lambda_\theta, \lambda_\beta)$ be the adjoint, and 
\begin{equation}
\begin{split}
& H(q, \lambda, u) = \lambda_x \dot x + \lambda_y \dot y + \lambda_\theta \dot \theta + \lambda_\beta \dot \beta \\
&= v(\lambda_x \cos\theta + \lambda_y \sin\theta - \lambda_\beta \sin\beta) + \omega(\lambda_\theta + \lambda_\beta),
\end{split}
\end{equation} 
the Hamiltonian. The Pontryagin Maximum Principle \cite{PonBolGamMis62} ensures that for every optimal trajectory-control pair $(q(t), u(t))$ over the time interval $[0, T]$, there exists an absolutely continuous vector-valued adjoint function $\lambda(t) \not= 0$ such that
\begin{align} \label{equ:lx}
\dot \lambda_x &= -\frac{\partial H}{\partial x} = 0, \\ \label{equ:ly}
\dot \lambda_y &= -\frac{\partial H}{\partial y} = 0, \\ \label{equ:ltheta}
\dot \lambda_\theta &= -\frac{\partial H}{\partial \theta} = v(\lambda_x \sin\theta - \lambda_y \cos\theta), \\ \label{equ:lbeta}
\dot \lambda_\beta &= -\frac{\partial H}{\partial \beta} = v\lambda_\beta \cos\beta,
\end{align}
and $H$ is maximized as a function of $v, \omega$ along the optimal trajectory, i.e. for all $t \in [0, T]$, 
\begin{equation}
H(q(t), \lambda(t), u(t)) = \max_{z \in U}\ H(q(t), \lambda(t), z). 
\end{equation}
Direct integration with help from (\ref{equ:system}) yields 
\begin{align} \label{equ:lambdax}
\lambda_x &= c_1, \\ \label{equ:lambday}
\lambda_y &= c_2, \\ \label{equ:lambdatheta}
\lambda_\theta &= c_1 y - c_2 x + c_3, \\ \label{equ:lambdabeta}
\lambda_\beta &= c_4 \exp\left(\int_{0}^{t} v \cos\beta(\tau)\ d\tau\right),
\end{align}
in which $c_1, c_2, c_3, c_4 = \lambda_\beta(0)$ are constant. Further integration of $\lambda_\beta$ requires more analysis which will be presented in the next sections.

Since the Pontryagin Maximum Principle guarantees that the Hamiltonian is maximized, as a function of the controls, we can define switching functions
\begin{align}\label{equ:switching1}
\phi_v &= \lambda_x \cos\theta + \lambda_y \sin\theta - \lambda_\beta \sin\beta, \\
\label{equ:switching2}
\phi_\omega &= \lambda_\theta + \lambda_\beta,
\end{align}
such that $v\phi_v = |\phi_v|$ and $\omega\phi_\omega = |\phi_\omega|$. That is $v = \mbox{sgn}(\phi_v)$ whenever $\phi_v \not = 0$ and $\omega = \mbox{sgn}(\phi_\omega)$ whenever $\phi_\omega \not = 0$.

An optimal trajectory is bound to satisfy the Pontryagin Maximum Principle, but not every trajectory that satifies the Pontryagin Maximum Principle is necessarily optimal. A trajectory that satisfies the Pontryagin Maximum Principle is called an \emph{extremal}. Extremal trajectories can be divided into two categories: \emph{regular} and \emph{singular}.

\section{Regular Extremals}
\begin{Definition}
An extremal trajectory is called \emph{regular} if the times at which $\phi_v = 0$ or $\phi_\omega = 0$ have zero measure, i.e. $\phi_v \not = 0$ and $\phi_\omega \not = 0$ almost everywhere. We define a \emph{regular primitive} to be a subtrajectory of a regular extremal, along which $\phi_v \not = 0$ and $\phi_\omega \not = 0$.  
\end{Definition}

Depending on the signs of the switching functions, there are four types of regular primitives:
\begin{itemize}
\item Forward-Left: $v = 1, \omega = 1$,
\item Forward-Right: $v = 1, \omega = -1$,
\item Backward-Left: $v = -1, \omega = 1$,
\item Backward-Right: $v = -1, \omega = -1$.
\end{itemize}
Let $[t_0, t_1] \subset [0, T]$ be the time interval over which a regular primitive is defined. Since the controls $v$ and $\omega$ are constant and $\dot \beta = \omega - v\sin\beta$, we can replace $d\tau$ with $d\beta / (\omega-v\sin\beta)$ in
(\ref{equ:lambdabeta}) to obtain
\begin{equation}
\begin{aligned}
\lambda_\beta(t) &= \lambda_\beta(t_0) \exp\left(\int_{\beta(t_0)}^{\beta(t)} \frac{v \cos\beta}{\omega-v\sin\beta}\ d\beta\right) \\
&= \lambda_\beta(t_0) \frac{\omega-v\sin\beta(t_0)}{\omega-v\sin\beta(t)}.
\end{aligned}
\end{equation} 
Note that $\beta$ will never reach $(v/\omega) \pi/2$ along a regular primitive if it does not start at $(v/\omega) \pi/2$, in which case our change of variables is valid. If $\beta(t_0) = (v/\omega) \pi/2$, then $\beta(t_0) \equiv (v/\omega) \pi/2$ and $\lambda_\beta(t) \equiv \lambda_\beta(t_0)$. Therefore along a regular primitive,
\begin{align}
\theta(t) &= \omega t + \theta(t_0), \\
x(t) &= x(t_0) + (v/\omega)(\sin(\theta) - \sin\theta(t_0)), \\
y(t) &= y(t_0) - (v/\omega)(\cos(\theta) - \cos \theta(t_0)), \\
\beta(t) &= 2 (v/\omega)\arctan \left(\frac{t-2v+K_1}{t+K_1}\right), \\
\lambda_\beta(t) &= \frac{K_2}{\omega-v\sin\beta(t)},
\end{align}
where constants $K_1$ and $K_2$ depend on the initial configuration by 
\begin{eqnarray}
K_1 &=& \frac{2}{v - \omega\tan (\beta(t_0)/2)}, \\
K_2 &=& \lambda_\beta(t_0) (\omega - v\sin\beta(t_0)).
\end{eqnarray}

\section{Singular Extremals}
\begin{Definition}
Extremals along which both $\phi_v \equiv 0$ and $\phi_\omega \equiv 0$ are called \emph{abnormal}. An extremal trajectory is called \emph{singular} if it contains a positive-measure subtrajectory along which $\phi_v \equiv 0$ or $\phi_\omega \equiv 0$ but not both. We call such a subtrajectory a \emph{singular primitive}. 
\end{Definition}
Chyba and Sekhavat characterized singular primitives \cite{ChySek99}. In particular, they showed that abnormal extremals are either trivial or satisfy $v \equiv 0, \omega = \pm 1$, which coincides with a special case of $\phi_v$-singular primitives below. We will summarize Chyba and Sekhavat's results below with a slight correction for $\phi_\omega$-singular primitives. 

\paragraph*{$\phi_v$-singular primitives.} $\phi_v \equiv 0$ on $[t_0, t_1]$: in this case, $c_1 = c_2 = c_4 = 0$ and $\lambda = (0, 0, c_3, 0)$. Controls are $\omega = \pm 1$ and $v \in [-1, 1]$ arbitrarily.

\paragraph*{$\phi_\omega$-singular primitives.} $\phi_\omega \equiv 0$ on $[t_0, t_1]$: in this case, $\phi_v = cte \not = 0$ and $\lambda_\beta = -\lambda_\theta$. 
Since $\phi_\omega \equiv 0$,
\begin{equation}
\begin{aligned}
\dot \phi_\omega &= \dot \lambda_\theta + \dot \lambda_\beta \\
&= v(c_1\sin\theta - c_2\cos\theta + \lambda_\beta \cos\beta) \equiv 0.
\end{aligned}
\end{equation}
Since $v = \mbox{sgn}(\phi_v) \not = 0$, $c_1\sin\theta - c_2\cos\theta + \lambda_\beta \cos\beta \equiv 0$. Taking derivative with respect to time, we obtain
\begin{equation}
\begin{aligned}
0 &\equiv (c_1\cos\theta + c_2\sin\theta)\dot \theta + \dot \lambda_\beta \cos\beta - \lambda_\beta \dot \beta \sin\beta\\
&= \omega(c_1\cos\theta + c_2\sin\theta -\lambda_\beta \sin\beta) + v\lambda_\beta\\
&= \omega\phi_v + v\lambda_\beta\\
&= \omega|\phi_v| + \lambda_\beta,
\end{aligned}
\end{equation}
using the fact that $v = \mbox{sgn}(\phi_v)$ and (\ref{equ:system}), (\ref{equ:lbeta}), (\ref{equ:lambdax}), (\ref{equ:lambday}), (\ref{equ:switching1}), and (\ref{equ:switching2}). Hence,
\begin{equation}\label{equ:singcontrol}
\begin{aligned}
\omega &= -\frac{\lambda_\beta}{|\phi_v|}\\
&= \frac{\lambda_\theta}{|\phi_v|}\\
&= \frac{c_1 y - c_2 x + c_3}{|\phi_v|},
\end{aligned}
\end{equation}
using the fact that $\lambda_\beta = -\lambda_\theta$ and (\ref{equ:lambdatheta}). Note that Chyba and Sekhavat have made a minus sign mistake in these calculations \cite{ChySek99}. Next, we analytically integrate the system over a $\phi_\omega$-singular primitive.

It is obvious from (\ref{equ:system}) that
\begin{align}
\dot x &= \omega \frac{d x}{d \theta} = v \cos\theta,\\
\dot y &= \omega \frac{d y}{d \theta} = v \sin\theta.
\end{align}
Hence using (\ref{equ:system}) and (\ref{equ:singcontrol}), we obtain
\begin{align}
(c_1 y - c_2 x + c_3)\frac{d x}{d \theta} &= \phi_v \cos\theta,\\
(c_1 y - c_2 x + c_3)\frac{d y}{d \theta} &= \phi_v \sin\theta.
\end{align}
Integrating with respect to $\theta$, we get
\begin{eqnarray} \label{equ:id1}
\begin{aligned}
-\frac{1}{2}c_2 x^2|_{t_0}^t & + c_3 x|_{t_0}^t + c_1 \int_{\theta(t_0)}^{\theta(t)} y \frac{d x}{d \theta} d \theta \\
&= \phi_v \left\{\sin\theta(t) - \sin\theta(t_0)\right\},
\end{aligned} \\ \label{equ:id2}
\begin{aligned}
 \frac{1}{2}c_1 y^2|_{t_0}^t & + c_3 y|_{t_0}^t - c_2 \int_{\theta(t_0)}^{\theta(t)} x\frac{d y}{d \theta} d \theta \\
&= \phi_v \left\{\cos\theta(t_0) - \cos\theta(t)\right\}.
\end{aligned}
\end{eqnarray}
Multiply (\ref{equ:id1}) by $c_2$ and add to (\ref{equ:id2}) multiplied by $-c_1$ to obtain
\begin{equation}\label{equ:bigid}
\begin{aligned}
& -\frac{1}{2}c^2_2 x^2|_{t_0}^t + c_2 c_3 x|_{t_0}^t - \frac{1}{2}c^2_1 y^2|_{t_0}^t - c_1 c_3 y|_{t_0}^t \\
& + c_1 c_2 \int_{\theta(t_0)}^{\theta(t)} \frac{d (xy)}{d \theta} d \theta \\
&= -\frac{1}{2}c^2_2 x^2|_{t_0}^t + c_2 c_3 x|_{t_0}^t - \frac{1}{2}c^2_1 y^2|_{t_0}^t - c_1 c_3 y|_{t_0}^t \\
& + c_1 c_2 (xy)|_{t_0}^t \\
&= -\frac{1}{2}(c_1y - c_2 x)^2|_{t_0}^t  - c_3(c_1 y - c_2 x)|_{t_0}^t \\
&= -\frac{1}{2}(\lambda^2_\theta - c_3^2)|_{t_0}^t \\
&= \frac{1}{2}\left\{\lambda^2_\theta(t_0) - \lambda^2_\theta(t)\right\} \\
&= \phi_v \left\{c_1 \cos\theta(t) + c_2\sin\theta(t) - c_1 \cos\theta(t_0) - c_2\sin\theta(t_0)\right\}\\
&= \phi_v \left\{\lambda_\theta(t_0)\sin \beta(t_0) - \lambda_\theta(t)\sin \beta(t)\right\}.
\end{aligned}
\end{equation}
Therefore, the following lemma holds.
\begin{lemma}
Consider a $\phi_\omega$-singular primitive over the time interval $[t_0, t_1]$. In that case for $t\in[t_0, t_1]$,
\begin{equation}
\begin{aligned}
\lambda_\theta(t) &= \phi_v\sin\beta(t) \\
&\pm \sqrt{\lambda^2_\theta(t_0) - 2 \phi_v\lambda_\theta(t_0)\sin\beta(t_0) + \phi^2_v\sin^2\beta(t)}.
\end{aligned}
\end{equation}
\end{lemma}

\begin{figure}[t!]
\begin{center}
\includegraphics[scale=0.52]{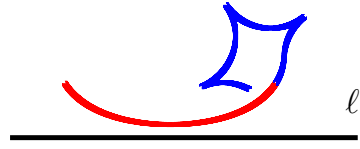}
\end{center}
\caption{A singular extremal trajectory composed of a $\phi_\omega$-singular primitive, plotted in \textcolor{red}{red}, and regular primitives, plotted in \textcolor{blue}{blue}. In this case, the $\phi_\omega$-singular primitive does not contain a straight line segment. The line $\ell$ is defined by $c_1 y - c_2 x + c_3 = 0$. } \label{fig:sing}
\end{figure}

Figure \ref{fig:sing} illustrates a singular extremal containing a $\phi_\omega$-singular primitive. Particularly, we are interested in those extremals that contain straight line motions, i.e. $v=\pm 1, \omega = 0$, since they are the only ones that can optimally take the robot to far destinations. A straight motion occurs only when $\lambda_\theta \equiv 0$ and $\theta \equiv \arctan(c_2/c_1)$ or $\theta \equiv \pi+\arctan(c_2/c_1)$. For the trajectory to cross the line $\lambda_\theta = 0$ at some time in $[t_0, t_1]$, (\ref{equ:bigid}) necessitates that either $\lambda_\theta(t_0) = 0$ or
\begin{equation}\label{equ:singphi}
\phi_v = \frac{\lambda_\theta(t_0)}{2\sin\beta(t_0)}.
\end{equation}
This is necessary but not sufficient as the control has to be within the valid range, i.e. $|\omega| \leq 1$, for the trajectory to remain singular and reach the line $\lambda_\theta = 0$.
In particular, if $\lambda_\theta(t_0) \not = 0$, then $\sin\beta(t_0) \not = 0$ and for $|\omega| \leq 1$ we must have $|\sin\beta(t_0)| \leq \frac{1}{2}$ which is equivalent to
\begin{equation}\label{equ:betarange}
-\frac{\pi}{6} \leq \beta(t_0) \leq \frac{\pi}{6} \mbox{ or } \frac{5\pi}{6} \leq \beta(t_0) \leq \frac{7\pi}{6}.
\end{equation}

\begin{figure*}[t!]
\begin{center}
\includegraphics[scale=0.52]{singular4.ps}
\end{center}
\caption{A singular extremal trajectory composed of a $\phi_\omega$-singular primitive, plotted in \textcolor{red}{red}, and regular primitives, plotted in \textcolor{blue}{blue}. In this case, the $\phi_\omega$-singular primitive contains a straight line segment. The line $\ell$ is defined by $c_1 y - c_2 x + c_3 = 0$. } \label{fig:singstraight2}
\end{figure*}

\begin{figure}[h!]
\begin{center}
\includegraphics[scale=0.52]{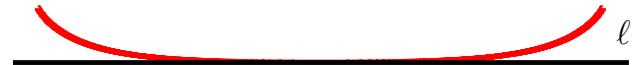}
\end{center}
\caption{A \emph{merging curve} followed by a straight motion and another merging curve. This case is a $\phi_\omega$-singular primitive that contains a segment of the line $\ell:\ c_1 y - c_2 x + c_3 = 0$.} \label{fig:singstraight1}
\end{figure}

For the moment assume $\lambda_\theta(t_0) \not = 0$.
Note that from (\ref{equ:switching1}) we have
\begin{equation}\label{equ:phiv}
\phi_v = c_1 \cos\theta + c_2 \sin\theta + \lambda_\theta \sin \beta,
\end{equation} 
which yields
\begin{equation}
\phi_v = c_1 \cos\theta(t_0) + c_2 \sin\theta(t_0) + \lambda_\theta(t_0) \sin \beta(t_0).
\end{equation} 
Together with (\ref{equ:singphi}), this equation yields an equation for $\lambda_\theta(t_0)$ in terms of $\theta(t_0)$ and $\beta(t_0)$:
\begin{equation}
\lambda_\theta(t_0)\cos2\beta(t_0) = 2 \sin\beta(t_0)\left\{c_1 \cos\theta(t_0) + c_2 \sin\theta(t_0)\right\}.
\end{equation} 
Therefore,
\begin{equation}
\lambda_\theta(t_0) = \dfrac{2 \sin\beta(t_0)\{c_1 \cos\theta(t_0) + c_2 \sin\theta(t_0)\}}{\cos2\beta(t_0)}, 
\end{equation}
as (\ref{equ:betarange}) guarantees that $\cos2\beta(t_0) \not = 0$ and
\begin{equation} \label{equ:phiinit}
\phi_v = \dfrac{c_1 \cos\theta(t_0) + c_2 \sin\theta(t_0)}{\cos2\beta(t_0)},
\end{equation} 
from (\ref{equ:singphi}).

\begin{figure*}
\begin{center}
\includegraphics[scale=0.52]{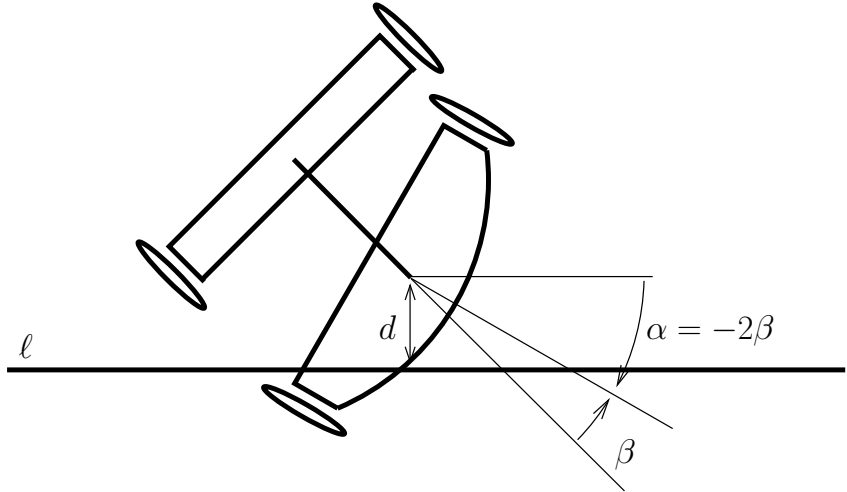}
\end{center}
\caption{The robot configuration along a merging curve. In this case, the robot tangentially joins the line $\ell:\ c_1 y - c_2 x + c_3 = 0$, $\alpha = \pm 2\beta$, $d = \pm 2 \sin \beta$, and $-\frac{\pi}{6} \leq \beta \leq \frac{\pi}{6}$ or $\frac{5\pi}{6} \leq \beta \leq \frac{7\pi}{6}$ along the trajectory.} \label{fig:distance}
\end{figure*}

Now for the trajectory to join the line $\lambda_\theta = 0$ tangentially so that it is able to follow a straight motion, such as in Figures \ref{fig:singstraight1} and \ref{fig:singstraight2}, the orientation vector $(\cos\theta, \sin\theta)$ at the tangency point should be $\pm (c_1, c_2)/\|(c_1, c_2)\|$. Hence, (\ref{equ:bigid}) yields
\begin{equation}
\phi_v \left\{\pm \|(c_1, c_2)\| - c_1 \cos\theta(t_0) - c_2\sin\theta(t_0)\right\} = \frac{1}{2}\lambda^2_\theta(t_0),
\end{equation}
from which we can derive two equations for $c_1, c_2,$ and $c_3$ in terms of $x(t_0), y(t_0), \theta(t_0),$ and $\beta(t_0)$:
\begin{equation}
\begin{split}
\phi_v & \left\{\pm \|(c_1, c_2)\| - c_1 \cos\theta(t_0) - c_2\sin\theta(t_0)\right\} = \frac{1}{2}\lambda^2_\theta(t_0) \\
&= \frac{1}{2}\left\{c_1 y(t_0) - c_2 x(t_0) + c_3\right\}^2,
\end{split}
\end{equation}
and from (\ref{equ:phiv}) and (\ref{equ:phiinit}),
\begin{equation}
\phi_v = \pm \| (c_1, c_2)\| = \dfrac{\| (c_1, c_2)\| \cos\alpha}{\cos2\beta(t_0)},
\end{equation}
in which $\alpha$ is the angle between the orientation vector at $t_0$, i.e. $(\cos\theta(t_0), \sin\theta(t_0))$, and $(c_1, c_2)$. Figure \ref{fig:distance} illustrates $\alpha$. It is clear then that
\begin{equation}
\cos2\beta(t_0) = \pm \cos\alpha,
\end{equation}
which implies
\begin{equation}
\alpha = \pm 2\beta(t_0) \mbox{ or } \pi \pm 2\beta(t_0).
\end{equation}
Denote the line $\lambda_\theta = 0$ by $\ell$, and the signed distance from $(x, y)$ to $\ell$ by 
\begin{equation}
d = \frac{c_1 y - c_2 x + c_3}{\|(c_1, c_2)\|}. 
\end{equation}
Figure \ref{fig:distance} illustrates $d$. From (\ref{equ:lambdatheta}) and (\ref{equ:singphi}), we obtain
\begin{equation}
|d(t_0)| = 2|\sin\beta(t_0)|.
\end{equation}
Moreover, the extremal angular control in (\ref{equ:singcontrol}) becomes
\begin{equation}
\omega = \frac{c_1 y - c_2 x + c_3}{\|(c_1, c_2)\|} = d.
\end{equation}
It is now obvious that for $|\omega| \leq 1$, the distance of the robot from $\ell$ has to remain within $[-1, 1]$. Since $t_0$ was not assumed to be the time at which a singular primitive starts, our results above hold for all $t \in [t_0, t_1]$ and hence the following lemma.
\begin{lemma}\label{lem:merging}
Let $\ell$ denote the line $c_1 y - c_2 x + c_3 = 0$ in the plane, and consider a $\phi_\omega$-singular primitive over the time interval $[t_0, t_1]$ that contains a straight line segment. First, all of the straight line segments of the considered primitive lie on $\ell$. Second, either $c_4 = 0$ and the entire primitive is a straight line, i.e.
\begin{equation}
\begin{aligned}
\alpha \equiv 0 \mbox{ or } \pi, \\
\omega \equiv 0, \\
\lambda_\beta \equiv \lambda_\theta \equiv 0, 
\end{aligned}
\end{equation}
or
\begin{align}
\alpha(t) &= \pm 2\beta(t), \\
\omega(t) &= d = \pm 2\sin\beta(t),\\
-\frac{\pi}{6} \leq \beta(t) \leq &\frac{\pi}{6} \mbox{ or } \frac{5\pi}{6} \leq \beta(t) \leq \frac{7\pi}{6},
\end{align}
in which $d$ is the signed distance from the center of robot to $\ell$, and $\alpha$ is the angle between the robot orientation and $\ell$. Figure \ref{fig:distance} illustrates the latter.
\end{lemma}
The former case above corresponds to $(\lambda_x, \lambda_y, \lambda_\theta, 0)$ adjoints which pertain to the Reeds-Shepp extremals \cite{SusTan91} including the Reeds-Shepp curves \cite{ReeShe90}. It is clear that the Reeds-Shepp curves are time optimal for our system. The latter case corrsponds to planar elastica connecting turning to straight segments \cite{Jur92}. 
\begin{Definition}
We call the latter case in Lemma \ref{lem:merging} a \emph{merging curve}.
\end{Definition}

\section{Conclusion and Future Work}

Although a complete characterization of optimal trajectories and control synthesis remains for future, we paved the way for such a complete characterization in this work. Based on this work, it is easy now to prove that for far enough destinations the optimal trajectory has to contain a straight line motion. Particularly, regular primitives and $\phi_v$-singular extremals cannot optimally take the robot to far destinations as they impose constant rotation ($\omega = \pm 1$) which becomes non-optimal after a certain period. We characterized the $\phi_\omega$-singular extremals which are the only extremals that can possibly contain a straight line segment. We gave a detailed characterization of \emph{merging curves} that connect regular extremals to straight line segments. It is obvious that the Reeds-Shepp curves are optimal for our system if the initial and goal $\beta$ are unimportant. The only optimal trajectories that contain a straight line segment and do not use merging curves are the Reeds-Shepp extremals.

\bibliographystyle{plain}
\bibliography{pub,main}

\end{document}